# A Generalization of Classical Symmetric Orthogonal Functions Using a Symmetric Generalization of Sturm-Liouville Problems

**Mohammad Masjed-Jamei**

*Department of Mathematics, K.N.Toosi University of Technology, P.O.Box 16315-1618, Tehran, Iran,
E-mail: mmjamei@kntu.ac.ir , mmjamei@yahoo.com*

**Abstract.** In this paper, usual Sturm-Liouville problems are extended for symmetric functions so that the corresponding solutions preserve the orthogonality property. Two basic examples, which are special cases of a generalized Sturm-Liouville problem, are then introduced. First example generalizes the associated Legendre functions having extensive applications in physics and the second example introduces a generic differential equation with various sub-cases that have orthogonal solutions. For instance, this generic equation possesses a symmetric differential equation containing a basic solution of symmetric orthogonal polynomials.

**1. Introduction**

As we know, partial differential equations can be reduced to ordinary equations by the method of separation of variables. For instance, this method is applied for solving partial differential equations of the general form

$$\rho(x, y, z)\left( A^*(t)\frac{\partial^2 u}{\partial t^2} + B^*(t)\frac{\partial u}{\partial t} \right) = Lu, \qquad (1.1)$$

where

$$Lu = \text{div}\,(k(x, y, z)\,\text{grad}\,u) - q(x, y, z)u \,.$$

This equation describes the propagation of a vibration such as electromagnetic or acoustic waves if $A^*(t) = 1$ and $B^*(t) = 0$. Moreover, when $A^*(t) = 0$ and $B^*(t) = 1$, it describes transfer processes, such as heat transfer or the diffusion of particles in a medium; and finally if $A^*(t) = 0$ and $B^*(t) = 0$, it describes the corresponding time-independent processes.

To have a unique solution for equation (1.1), which corresponds to an actual physical problem, it is clear that we must impose some supplementary conditions. The most typical conditions are initial or boundary conditions. The initial conditions for equation (1.1) are usually the values of $u(x,y,z,t)$ and $\partial u(x,y,z,t)/\partial t$, while the simplest boundary conditions have the form

$$\left[\alpha(x,y,z)u + \beta(x,y,z)\frac{\partial u}{\partial \eta}\right]_S = 0, \qquad (1.2)$$

where $\alpha(x,y,z)$ and $\beta(x,y,z)$ are given functions; $S$ is the surface bounding the domain where (1.1) is to be solved; $\partial u/\partial \eta$ is the derivative in the direction of the outward normal to $S$. Hence, particular solutions of (1.1) under the boundary conditions (1.2) will be found if one looks for a solution of the form

$$u(x,y,z,t) = T(t)v(x,y,z). \qquad (1.3)$$

By substituting (1.3) into the main equation (1.1) one respectively gets

$$A^*(t)T'' + B^*(t)T' + \lambda T = 0, \qquad (1.4)$$
$$Lv + \lambda \rho v = 0, \qquad (1.5)$$

where $\lambda$ is a constant. Clearly equation (1.4) can be solved for typical problems in mathematical physics. However, to solve (1.5) we are to use a boundary condition that follows from (1.2), namely $(\alpha(x,y,z)v + \beta(x,y,z)\frac{\partial v}{\partial \eta})\big|_S = 0$. This problem is in fact a multidimensional boundary value problem. On the other hand, a multidimensional problem can be simplified to a one-dimensional problem if (1.5) is reduced to an equation of the form

$$Ly + \lambda \rho(x) y = 0, \qquad (1.6)$$

where

$$Ly = \frac{d}{dx}\left(k(x)\frac{dy}{dx}\right) - q(x)y, \quad k(x) > 0, \quad \rho(x) > 0. \qquad (1.6.1)$$

The equation (1.6) should be considered on an open interval, say $(a,b)$, with boundary conditions in the form

$$\begin{aligned}\alpha_1 y(a) + \beta_1 y'(a) &= 0, \\ \alpha_2 y(b) + \beta_2 y'(b) &= 0,\end{aligned} \qquad (1.7)$$

where $\alpha_1, \alpha_2$ and $\beta_1, \beta_2$ are given constants and $k(x), k'(x), q(x)$ and $\rho(x)$ in (1.6) are to be assumed continuous for $x \in [a,b]$. The simplified boundary value problem (1.6)-(1.7) is called a regular Sturm-Lioville problem. Moreover, if one of the boundary points $a$ and $b$ is singular (i.e. $k(a)=0$ or $k(b)=0$), the problem is transformed to a singular Sturm-Liouville problem. In this case, we can ignore boundary conditions (1.7) and obtain the orthogonality property directly.

Now, let $y_n(x)$ and $y_m(x)$ be two eigenfunctions of equation (1.6). According to Sturm-Liouville theory [3,9], these functions are orthogonal with respect to the weight function $\rho(x)$ on $(a,b)$ under the given conditions (1.7), i.e.

$$\int_a^b \rho(x) y_n(x) y_m(x)\, dx = \left( \int_a^b \rho(x) y_n^2(x)\, dx \right) \delta_{n,m} \quad \text{if} \quad \delta_{n,m} = \begin{cases} 0 & if \quad n \neq m \\ 1 & if \quad n = m \end{cases}. \tag{1.8}$$

Many important special functions in theoretical and mathematical physics are the solutions of regular or singular Sturm-Liouville problems that satisfy the orthogonality condition (1.8). For instance, the associated Legendre functions [3], Bessel functions [3,9], Fourier trigonometric sequences [1], Ultraspherical functions [5,10], Hermite functions [5,10] and so on are particular solutions of Sturm-Liouville problems. Fortunately, most of these mentioned functions have the symmetry property, namely $\Phi_n(-x) = (-1)^n \Phi_n(x)$, and because of this property they have found various applications in physics, see e.g. [3,9] for more details. Now, if we can extend the mentioned examples symmetrically and preserve their orthogonality property, it seems that we will be able to find some new applications in physics that logically extend the previous applications. In this research, by extending the Sturm-Liouville problems for symmetric functions, we generalize some classical symmetric orthogonal functions and obtain their orthogonality property. Hence, the construction of the paper would be as follows: In section 2, we introduce some initial conditions in order to be able to extend the Sturm-Lioville problem (1.6) with symmetric solutions. In this way, a main theorem is given. After proving the main theorem, two practical examples of generalized Sturm-Liouville problems type, i.e. a symmetric generalization of the associated Legendre functions and a generic differential equation with an explicit solution of symmetric orthogonal functions including almost all known symmetric orthogonal polynomials, such as Ultraspherical polynomials, generalized Ultraspherical polynomials (GUP), Hermite polynomials, generalized Hermite polynomials (GHP) and two further sequences of *finite* symmetric orthogonal polynomials are introduced in sections 3 and 4. Here we should add that finite orthogonal polynomials are a less-known subject in the theory of special functions. So, to study this topic we refer the readers to e.g. [7].

## 2. Generalized Sturm-Liouville problems with symmetric solutions

Without loss of generality, let $y = \Phi_n(x)$ be a sequence of symmetric functions satisfying the following differential equation

$$A(x)\Phi_n''(x) + B(x)\Phi_n'(x) + \left( \lambda_n C(x) + D(x) + \frac{1-(-1)^n}{2} E(x) \right) \Phi_n(x) = 0, \tag{2.1}$$

where $A(x)$, $B(x)$, $C(x)$, $D(x)$ and $E(x)$ are independent functions and $\{\lambda_n\}$ is a sequence of constants. Clearly choosing $E(x) = 0$ in equation (2.1) is equivalent to the same Sturm-Liouville equation as (1.6). Since $\Phi_n(x)$ is a symmetric sequence, by replacing the symmetry property $\Phi_n(-x) = (-1)^n \Phi_n(x)$ into equation (2.1), one can immediately conclude that $A(x)$, $C(x)$, $D(x)$ and $E(x)$ are *even* functions, while $B(x)$ must be an *odd* function. This note will be used in the paper frequently.

To prove the orthogonality property of the sequence $\Phi_n(x)$, both sides of equation (2.1) should first be multiplied by the positive function

$$R(x) = \exp(\int \frac{B(x) - A'(x)}{A(x)} dx) = \frac{1}{A(x)} \exp(\int \frac{B(x)}{A(x)} dx), \quad (2.2)$$

in order to convert it in the form of a self-adjoint differential equation. Note that $R(x)$ is an even function in this relation, because $B(x)$ is odd and $A(x)$ is even. Therefore, the self-adjoint form of equation (2.1) becomes

$$\frac{d}{dx}(A(x) R(x) \frac{d\Phi_n(x)}{dx}) = -(\lambda_n C(x) + D(x)) R(x) \Phi_n(x) - \frac{1 - (-1)^n}{2} E(x) R(x) \Phi_n(x). \quad (2.3)$$

Since $A(x) R(x)$ is an even function, the orthogonality interval corresponding to equation (2.3) should be considered symmetric, say $[-\theta, \theta]$. Hence, by assuming that $x = \theta$ is a root of the function $A(x) R(x)$ and applying the Sturm-Liouville theorem on equation (2.3) we have

$$[A(x) R(x) (\Phi'_n(x) \Phi_m(x) - \Phi'_m(x) \Phi_n(x))]_{-\theta}^{\theta} =$$
$$(\lambda_m - \lambda_n) \int_{-\theta}^{\theta} C(x) R(x) \Phi_n(x) \Phi_m(x) dx - (\frac{(-1)^m - (-1)^n}{2}) \int_{-\theta}^{\theta} E(x) R(x) \Phi_n(x) \Phi_m(x) dx. \quad (2.4)$$

Obviously the left-hand side of (2.4) is zero. So, to prove the orthogonality property, it remains to show that the following value

$$F(n, m) = \frac{(-1)^m - (-1)^n}{2} \int_{-\theta}^{\theta} E(x) R(x) \Phi_n(x) \Phi_m(x) dx, \quad (2.5)$$

is always equal to zero for every $m, n \in Z^+$. For this purpose, generally four cases should be considered for *m* and *n*:

a) If both *m* and *n* are even (or odd), then it is clear that $F(n,m) = 0$, because we have $F(2i, 2j) = F(2i+1, 2j+1) = 0$.

b) If one of the two mentioned values is odd and the other one is even (or conversely) then

$$F(2i, 2j+1) = -\int_{-\theta}^{\theta} E(x) R(x) \Phi_{2i}(x) \Phi_{2j+1}(x) dx. \tag{2.6}$$

But $E(x)$, $R(x)$ and $\Phi_{2i}(x)$ are even functions and $\Phi_{2j+1}(x)$ is odd in above relation. So the integrand of (2.6) is an odd function and consequently $F(2i, 2j+1) = 0$. This issue also holds for the case $n = 2i+1$ and $m = 2j$, i.e. $F(2i+1, 2j) = 0$.

By noting these comments, the main theorem of the paper can here be expressed.

**Theorem 1.** The symmetric sequence $\Phi_n(x) = (-1)^n \Phi_n(-x)$, as a specific solution of differential equation (2.1), satisfies the orthogonality relation

$$\int_{-\theta}^{\theta} W^*(x) \Phi_n(x) \Phi_m(x) dx = \left( \int_{-\theta}^{\theta} W^*(x) \Phi_n^2(x) dx \right) \delta_{n,m}, \tag{2.7}$$

where

$$W^*(x) = C(x) R(x) = C(x) \exp\left( \int \frac{B(x) - A'(x)}{A(x)} dx \right), \tag{2.8}$$

denotes the corresponding weight function. Note that $W^*(x)$ in (2.8) must be a positive and even function on $[-\theta, \theta]$.

Now, it is a good position to present some practical examples that generalize the well-known symmetric orthogonal functions and preserve their orthogonality property. The next sections 3 and 4 introduce two basic examples in this regard.

## 3. A symmetric generalization of the associated Legendre functions.

When Laplace's equation $\nabla^2 u = 0$ (or Helmholtz's equation $\nabla^2 u = \lambda u$) is solved in spherical coordinates $r, \theta, \varphi$, the following results are derived [9]

$$\nabla^2 u = \Delta_r u + \frac{1}{r^2} \Delta_{\theta,\varphi} u, \tag{3.1}$$

where

$$\Delta_r u = \frac{1}{r^2} \frac{\partial}{\partial r} \left( r^2 \frac{\partial u}{\partial r} \right),$$

$$\Delta_{\theta,\varphi} u = \frac{1}{\sin \theta} \frac{\partial}{\partial \theta} \left( \sin \theta \frac{\partial u}{\partial \theta} \right) + \frac{1}{\sin^2 \theta} \frac{\partial^2 u}{\partial \varphi^2}. \tag{3.2}$$

As we said in the introduction, by separating $u = R(r) \Psi(\theta) \Omega(\varphi)$ and substituting it into Laplace's equation (Potential equation) one gets three ordinary differential equations

$$(r^2 R')' = \mu\ R,$$
$$\Omega'' + \upsilon\ \Omega = 0, \tag{3.3}$$
$$\sin \theta \frac{d}{d\theta} \left( \sin \theta \frac{d\Psi}{d\theta} \right) + (\mu \sin^2 \theta - \upsilon)\ \Psi = 0,$$

where $\mu$ and $\upsilon$ are two constant values. For $x = \cos\theta$, the last equation of (3.3) changes to

$$(1-x^2)\frac{d^2\Psi}{dx^2} - 2x\frac{d\Psi}{dx} + (\mu - \frac{\upsilon}{1-x^2})\Psi = 0 , \qquad (3.4)$$

which is called the associated Legendre differential equation [3] and its solutions are consequently known as the associated Legendre functions. Hence, these functions play a key role in the potential theory. Generally there are three types of solution for equation (3.4) depending on different values of $\mu$ and $\upsilon$:

i) If $\mu = (n+\alpha)(n+\alpha+1)$ and $\upsilon = \alpha^2$; $n \in Z^+$, $\alpha > -1$, then the related solution is indicated by

$$\Psi(x) = U_n^{(\alpha)}(x) = (1-x^2)^{\frac{\alpha}{2}} P_n^{(\alpha,\alpha)}(x), \qquad (3.5)$$

where $P_n^{(\alpha,\alpha)}(x)$, known as Ultraspherical polynomials, is a special case of Jacobi polynomials [5,10]

$$P_n^{(\alpha,\beta)}(x) = \sum_{k=0}^{n} \binom{n+\alpha+\beta+k}{k}\binom{n+\alpha}{n-k}(\frac{x-1}{2})^k , \qquad (3.6)$$

for $\alpha = \beta$. In this way, the solution (3.5) satisfies the orthogonality relation

$$\int_{-1}^{1} U_n^{(\alpha)}(x) U_m^{(\alpha)}(x) dx = (\int_{-1}^{1} (U_n^{(\alpha)}(x))^2 dx)\delta_{n,m} = \frac{2^{2\alpha+1}\Gamma^2(n+\alpha+1)}{n!(2n+2\alpha+1)\Gamma(n+2\alpha+1)} \delta_{n,m}. \qquad (3.7)$$

ii) If $\mu = n(n+1)$ and $\upsilon = m^2$; $m,n \in Z^+$, (3.4) has a solution as

$$\Psi(x) = P_n^m(x) = (1-x^2)^{\frac{m}{2}} \frac{d^m(P_n(x))}{dx^m}, \qquad (3.8)$$

where $P_n(x) = P_n^{(0,0)}(x)$ denotes the Legendre polynomials [9]. Moreover, according to [3], (3.8) satisfies the relation

$$\int_{-1}^{1} P_i^m(x) P_j^m(x) dx = (\int_{-1}^{1} (P_i^m(x))^2 dx)\delta_{i,j} = \frac{2(i+m)!}{(2i+1)(i-m)!} \delta_{i,j}. \qquad (3.9)$$

iii) Finally if $\mu = n(n+1)$ and $\upsilon = \alpha^2$; $-1 < \alpha < 1$ are selected in (3.4), the related solution takes the form

$$\Psi(x) = V_n^{(\alpha)}(x) = (\frac{1-x}{1+x})^{\frac{\alpha}{2}} P_n^{(\alpha,-\alpha)}(x), \qquad (3.10)$$

that satisfies the orthogonality relation

$$\int_{-1}^{1} V_n^{(\alpha)}(x) V_m^{(\alpha)}(x)\,dx = (\int_{-1}^{1} (V_n^{(\alpha)}(x))^2\,dx)\delta_{n,m} = \frac{2\Gamma(n+1+\alpha)\Gamma(n+1-\alpha)}{(n!)^2 (2n+1)} \delta_{n,m}. \quad (3.11)$$

Now, to extend the associated Legendre functions, it is enough to choose

$$A(x) = 1 - x^2 = A(-x) \quad ; \quad B(x) = -2x = -B(-x),$$
$$C(x) = 1 = C(-x) \quad ; \quad D(x) = -\frac{\upsilon}{1-x^2} = D(-x), \quad (3.12)$$
$$\lambda_n = \mu_n \quad ; \quad E(x) = E(-x) \quad \text{Arbitrary},$$

from the main equation (2.1) and apply the theorem 1 to establish the orthogonality property on the same interval [-1,1]. To do this task, without loss of generality, we suppose that $\Phi_n(x) = Q_n(x;\upsilon, E(x))$ is a known solution of the following differential equation

$$(1-x^2)\Phi_n''(x) - 2x\Phi_n'(x) + (\mu_n - \frac{\upsilon}{1-x^2} + \frac{1-(-1)^n}{2} E(x))\Phi_n(x) = 0. \quad (3.13)$$

According to theorem 1, one can show that

$$\int_{-1}^{1} Q_n(x;\upsilon, E(x)) Q_m(x;\upsilon, E(x))\,dx = (\int_{-1}^{1} (Q_n(x;\upsilon, E(x)))^2\,dx)\delta_{n,m}, \quad (3.14)$$

Provided that the arbitrary function $E(x)$ is *even*. Evidently various options can be selected for $E(x)$, which are directly related to the orthogonal sequence $Q_n(x;\upsilon, E(x))$. For example, choosing the even function $E(x) = 0$ gives one of the three cases of usual associated Legendre functions, stated above. The special case $E(x) = E(-x) = -2/x^2$ will also be introduced in the section 4, part (4.a).

**Remark 1.** Notice since $Q_n(x;\upsilon, E(x))$ (as a known solution of equation (3.13)) is now an orthogonal sequence under the predetermined condition $E(-x) = E(x)$, it can also be applied for the expansion of functions. As it is known, many boundary value problems of mathematical physics can be solved by using the expansion of functions in terms of the eigenfunctions of a Sturm-Liouville problem. Hence, if we assume

$$f(x) = \sum_{n=0}^{\infty} q_n Q_n(x;\upsilon, E(x)), \quad (3.15)$$

then according to (3.14), the unknown coefficients $q_n$ can be shown by

$$q_n = \int_{-1}^{1} f(x) Q_n(x;\upsilon, E(x))\,dx \Big/ \int_{-1}^{1} (Q_n(x;\upsilon, E(x)))^2\,dx. \quad (3.16)$$

## 4. A generic differential equation and its orthogonal solutions

In this section, by using the main differential equation (2.1) and subsequently theorem 1, we are going to introduce a fundamental class of symmetric orthogonal functions that generates almost all known symmetric sequences of orthogonal polynomials. For this purpose, let us start with a second order differential equation, which is considered in [8], as (4.1)

$$x^2(px^2+q)\Phi_n''(x) + x(rx^2+s)\Phi_n'(x) - \left(n(r+(n-1)p)x^2 + (1-(-1)^n)s/2\right)\Phi_n(x) = 0.$$

Clearly this equation is a special case of (2.1) for

$$\begin{aligned}
A(x) &= x^2(px^2+q) & ; & \quad An \text{ even function}, \\
B(x) &= x(rx^2+s) & ; & \quad An \text{ odd function}, \\
C(x) &= x^2 & ; & \quad An \text{ even function}, \qquad (4.2)\\
D(x) &= 0 & ; & \quad An \text{ even function}, \\
E(x) &= -s & ; & \quad An \text{ even function},
\end{aligned}$$

and $\lambda_n = -n(r+(n-1)p)$. According to [8], the polynomial solution of equation (4.1) is (4.3)

$$\Phi_n(x) = S_n\!\left(\begin{array}{cc} r & s \\ p & q \end{array}\bigg| x\right) = \sum_{k=0}^{[n/2]} \binom{[n/2]}{k} \left(\prod_{i=0}^{[n/2]-(k+1)} \frac{(2i+(-1)^{n+1}+2[n/2])p+r}{(2i+(-1)^{n+1}+2)q+s}\right) x^{n-2k},$$

where neither both values $q$ and $s$ nor both values $p$ and $r$ can vanish together. From (4.3) it can be concluded that the leading coefficient of polynomials is

$$K_n = \prod_{i=0}^{[n/2]-1} \frac{(2i+(-1)^{n+1}+2\,[n/2])\,p+r}{(2i+(-1)^{n+1}+2)\,q+s}. \qquad (4.3.1)$$

Therefore, for instance, the monic polynomials of (4.3) for $n=0,\ldots,5$ take the forms

$$\begin{aligned}
&\bar{S}_0\!\left(\begin{array}{cc} r & s \\ p & q \end{array}\bigg| x\right) = 1 \qquad\qquad ; \qquad \bar{S}_1\!\left(\begin{array}{cc} r & s \\ p & q \end{array}\bigg| x\right) = x, \\
&\bar{S}_2\!\left(\begin{array}{cc} r & s \\ p & q \end{array}\bigg| x\right) = x^2 + \frac{q+s}{p+r} \quad ; \quad \bar{S}_3\!\left(\begin{array}{cc} r & s \\ p & q \end{array}\bigg| x\right) = x^3 + \frac{3q+s}{3p+r}x, \\
&\bar{S}_4\!\left(\begin{array}{cc} r & s \\ p & q \end{array}\bigg| x\right) = x^4 + 2\frac{3q+s}{5p+r}x^2 + \frac{(3q+s)(q+s)}{(5p+r)(3p+r)}, \\
&\bar{S}_5\!\left(\begin{array}{cc} r & s \\ p & q \end{array}\bigg| x\right) = x^5 + 2\frac{5q+s}{7p+r}x^3 + \frac{(5q+s)(3q+s)}{(7p+r)(5p+r)}x.
\end{aligned} \qquad (4.3.2)$$

Furthermore, in [8] we showed that the monic polynomials

$$\bar{S}_n\left(\begin{matrix}r & s \\ p & q\end{matrix}\bigg| x\right) = \frac{1}{K_n} S_n\left(\begin{matrix}r & s \\ p & q\end{matrix}\bigg| x\right), \qquad (4.3.3)$$

satisfy the recurrence relation

$$\bar{S}_{n+1}(x) = x\bar{S}_n(x) + C_n\left(\begin{matrix}r & s \\ p & q\end{matrix}\right)\bar{S}_{n-1}(x) \; ; \quad \bar{S}_0(x) = 1, \; \bar{S}_1(x) = x, \; n = 1,2,\dots, \qquad (4.4)$$

where

$$C_n\left(\begin{matrix}r & s \\ p & q\end{matrix}\right) = \frac{pq\, n^2 + \left((r-2p)q - (-1)^n ps\right)n + (r-2p)s(1-(-1)^n)/2}{(2pn+r-p)(2pn+r-3p)}. \qquad (4.4.1)$$

Since the form of recurrence relation (4.4) is completely available, Favard theorem [5] can here be applied to design a generic orthogonality relation for the monic polynomials (4.3.3) as follows (4.5)

$$\int_{-\theta}^{\theta} W\left(\begin{matrix}r & s \\ p & q\end{matrix}\bigg| x\right) \bar{S}_n\left(\begin{matrix}r & s \\ p & q\end{matrix}\bigg| x\right) \bar{S}_m\left(\begin{matrix}r & s \\ p & q\end{matrix}\bigg| x\right) dx = \left((-1)^n \prod_{i=1}^{n} C_i\left(\begin{matrix}r & s \\ p & q\end{matrix}\right) \int_{-\theta}^{\theta} W\left(\begin{matrix}r & s \\ p & q\end{matrix}\bigg| x\right) dx \right) \delta_{n,m},$$

where (4.5.1)

$$W\left(\begin{matrix}r & s \\ p & q\end{matrix}\bigg| x\right) = x^2 \exp(\int \frac{(r-4p)x^2 + (s-2q)}{x(px^2+q)} dx) = \exp(\int \frac{(r-2p)x^2 + s}{x(px^2+q)} dx),$$

is the related weight function. Note that according to (2.2) and (4.2), the function $A(x)R(x) = x^2(px^2+q)\exp(\int \frac{(r-2p)x^2+s}{x(px^2+q)} dx)$ must vanish at $x = \theta$ (and subsequently at $x = -\theta$). The standard values of this point are either $\theta = 1$ or $\theta = \infty$. Here let us add that the weight function (4.5.1), as an analogue of Pearson distributions family [9,10], can be investigated from statistical point of view. This function satisfies a first order differential equation of the form

$$x\frac{d}{dx}((px^2+q)W(x)) = (rx^2+s)W(x). \qquad (4.6)$$

Now, by noting the stated introduction in this section, we are in a good position to introduce a generic differential equation and its orthogonal solutions. Hence, without loss of generality, let us define the following equation

$$x^2(ax^\lambda + b)y'' + x(cx^\lambda + d)y' + \left(\alpha_n x^\lambda + \beta(1-(-1)^n)/2\right)y = 0, \qquad (4.7)$$

in which $a,b,c,d,\lambda$ and $\beta$ are real numbers and $\{\alpha_n\}$ is a sequence of constants. First it is clear that (4.7) is a special case of the main equation (2.1) for

$$A(x) = x^2(ax^\lambda + b) \quad ; \quad B(x) = x(cx^\lambda + d) \quad ;$$
$$C(x) = x^\lambda \quad\quad\quad\quad ; \quad D(x) = 0 \quad\quad\quad ; \quad E(x) = \beta. \tag{4.8}$$

But, our main aim is to find some conditions under which the general solution of equation (4.7) satisfies the conditions of theorem 1. For this purpose, we should reconsider the equation (4.1) and suppose $x = t^{\lambda/2}$. By substituting

$$y' = \frac{dy}{dx} = \frac{2}{\lambda} t^{-\frac{\lambda}{2}+1} \frac{dy}{dt},$$
$$y'' = \frac{d^2y}{dx^2} = \frac{4}{\lambda^2} t^{-\lambda+2} \frac{d^2y}{dt^2} + (\frac{4}{\lambda^2} - \frac{2}{\lambda}) t^{-\lambda+1} \frac{dy}{dt}, \tag{4.9}$$

into equation (4.1) and noting its general solution in (4.3), we get

(4.10)

$$t^2(pt^\lambda + q)\frac{d^2y}{dt^2} + t\left((\frac{\lambda}{2}r + (1-\frac{\lambda}{2})p)t^\lambda + \frac{\lambda}{2}s + (1-\frac{\lambda}{2})q\right)\frac{dy}{dt}$$
$$- \left(\frac{\lambda^2}{4}n(r + (n-1)p)t^\lambda + \frac{\lambda^2}{4}s\,(\frac{1-(-1)^n}{2})\right)y = 0 \Leftrightarrow y = S_n\begin{pmatrix} r & s \\ p & q \end{pmatrix}\bigg| t^{\frac{\lambda}{2}}\end{pmatrix}.$$

Here, for convenience, let us assume that

$$p = a \;,\; q = b \;,\; \frac{\lambda}{2}r + (1-\frac{\lambda}{2})p = c \;,\; \frac{\lambda}{2}s + (1-\frac{\lambda}{2})q = d. \tag{4.10.1}$$

Hence (4.10) is transformed to the same equation as (4.7) with the specific parameters

$$\alpha_n = -\frac{\lambda}{4} n(2c + (\lambda n - 2)a) \quad \text{and} \quad \beta = -\frac{\lambda}{4}(2d + (\lambda-2)b), \tag{4.11}$$

and the general solution (4.11.1)

$$y = S_n\begin{pmatrix} \frac{2}{\lambda}c + (1-\frac{2}{\lambda})a, & \frac{2}{\lambda}d + (1-\frac{2}{\lambda})b \\ a, & b \end{pmatrix}\bigg| x^{\frac{\lambda}{2}}\end{pmatrix}$$
$$= \sum_{k=0}^{[n/2]} \binom{[n/2]}{k} \left( \prod_{i=0}^{[n/2]-(k+1)} \frac{\left(2i + (-1)^{n+1} + 2\,[n/2] + 1 - 2/\lambda\right)a + 2c/\lambda}{\left(2i + (-1)^{n+1} + 3 - 2/\lambda\right)b + 2d/\lambda} \right) x^{(n-2k)\frac{\lambda}{2}}.$$

Now, the necessary conditions for the solution of equation corresponding to (4.11) in order to satisfy the conditions of theorem 1 are respectively that:

i) The general solution (4.11.1) must be symmetric.
ii) The coefficients of equation (4.7) must alternatively be even and odd.

Consequently we should have

$$\begin{cases}(-x)^{\frac{\lambda}{2}}=-x^{\frac{\lambda}{2}}\\(-x)^{\lambda}=x^{\lambda}\end{cases} \Rightarrow \begin{cases}(-1)^{\lambda/2}=-1\\(-1)^{\lambda}=+1\end{cases} \; ; \; \lambda \in R \, . \qquad (4.12)$$

Under these given conditions, many options can be selected for $\lambda$. For instance, if $\lambda = 2/3$ in (4.11), then the differential equation $\qquad\qquad\qquad\qquad\qquad (4.13)$

$$x^2(ax^{\frac{2}{3}}+b)y'' + x(cx^{\frac{2}{3}}+d)y' - \left(\frac{n}{3}(c+(\frac{n}{3}-1)a)x^{\frac{2}{3}} + \frac{1}{6}(d-\frac{2}{3}b)(1-(-1)^n)\right)y = 0,$$

has a general solution as

$$y = S_n\!\left(\begin{array}{cc}3c-2a, & 3d-2b\\a, & b\end{array}\bigg|\sqrt[3]{x}\right), \qquad (4.13.1)$$

which is symmetric and therefore satisfies the conditions of theorem 1. It also has a corresponding weight function as

$$W_1(x) = x^{\frac{2}{3}}\exp(\int \frac{(c-8a/3)x^{\frac{2}{3}}+(d-2b)}{ax^{5/3}+bx}dx). \qquad (4.13.2)$$

Finally, substituting $x = t^{1/3}$ into (4.5) gives the orthogonality relation of this sequence of orthogonal functions as $\qquad\qquad\qquad\qquad\qquad (4.13.3)$

$$\int_{-\theta^3}^{\theta^3} W_1(t)\bar{S}_n\!\left(\begin{array}{cc}3c-2a, & 3d-2b\\a, & b\end{array}\bigg|\sqrt[3]{t}\right)\bar{S}_m\!\left(\begin{array}{cc}3c-2a, & 3d-2b\\a, & b\end{array}\bigg|\sqrt[3]{t}\right)dt$$
$$= \left((-1)^n\prod_{i=1}^{n}C_i\!\left(\begin{array}{cc}3c-2a, & 3d-2b\\a, & b\end{array}\right)\int_{-\theta^3}^{\theta^3}W_1(t)dt\right)\delta_{n,m}\,.$$

Note that the right-hand side of (4.13.3) can be computed by referring to (4.4.1).
On the other hand, there are generally four special sub-classes of polynomials (4.3) (and consequently the main sequence (4.11.1) for $\lambda = 2$) that are orthogonal with respect to four special sub-solutions of equation (4.6). The main properties of these sub-classes are summarized in sections (4.a) to (4.d).

**4.a. First subclass: generalized Ultraspherical polynomials (GUP)**

These polynomials were first investigated by Chihara [5]. He obtained their main properties by employing a direct relation between them and Jacobi orthogonal polynomials. To study GUP individually we refer the reader to [2,4,6].
By replacing the initial values

$$(a,b,c,d,\lambda) = (-1,1,-2u-2v-2,2u,2), \qquad (4.14)$$

in the main sequence (4.11.1), the explicit form of monic GUP is derived as

(4.14.1)
$$\bar{S}_n\begin{pmatrix} -2u-2v-2, & 2u \\ -1, & 1 \end{pmatrix} x = \prod_{i=0}^{[n/2]-1} \frac{2i+2u+2-(-1)^n}{-2i-(2v+2u+2-(-1)^n+2[n/2])}$$
$$\times \sum_{k=0}^{[n/2]} \binom{[n/2]}{k} (\prod_{i=0}^{[n/2]-(k+1)} \frac{-2i-(2v+2u+2-(-1)^n+2[n/2])}{2i+2u+2-(-1)^n}) x^{n-2k}$$

that corresponds to the weight function

$$W\begin{pmatrix} -2u-2v-2 & 2u \\ -1 & 1 \end{pmatrix} x = x^{2u}(1-x^2)^v, \qquad (4.14.2)$$

and satisfy the orthogonality relation

$$\int_{-1}^{1} x^{2u}(1-x^2)^v \bar{S}_n\begin{pmatrix} -2u-2v-2, & 2u \\ -1, & 1 \end{pmatrix} x \bar{S}_m\begin{pmatrix} -2u-2v-2, & 2u \\ -1, & 1 \end{pmatrix} x dx$$
$$= \left((-1)^n \prod_{i=1}^{n} C_i\begin{pmatrix} -2u-2v-2, & 2u \\ -1, & 1 \end{pmatrix} \int_{-1}^{1} x^{2u}(1-x^2)^v dx\right) \delta_{n,m}, \qquad (4.14.3)$$

where

$$\int_{-1}^{1} x^{2u}(1-x^2)^v dx = B(u+\frac{1}{2}, v+1) = \frac{\Gamma(u+1/2)\Gamma(v+1)}{\Gamma(u+v+3/2)}. \qquad (4.14.4)$$

From (4.14.4), it can be deduced that the constraints on the parameters $u$ and $v$ are respectively $u+1/2>0$, $(-1)^{2u}=1$ and $v+1>0$. Note that $B(u,v)$ in this relation denotes the Beta integral [9] and $\Gamma(u)$ denotes the gamma function [3].

Finally the polynomials GUP satisfies the symmetric differential equation (4.14.5)

$$x^2(-x^2+1)\Phi_n''(x) - 2x((u+v+1)x^2-u)\Phi_n'(x) + \left(n(2u+2v+n+1)x^2 + ((-1)^n-1)u\right)\Phi_n(x) = 0$$

This equation is evidently a special case of the generic differential equation (4.7) for the initial vector (4.14) and the defined parameters (4.11).

But, as we reminded in section 3, there is a special case for $E(x)$ in equation (3.13) that can be derived from (4.4.15). To generate this special case, first let us define the following sequence

$$G_n(x;a,b) = x^a(1-x^2)^{\frac{b}{2}} \bar{S}_n\begin{pmatrix} -2a-2b-2, & 2a \\ -1, & 1 \end{pmatrix} x \quad ; \quad a > -\frac{1}{2}, \quad b > -1. \qquad (4.15)$$

If the above sequence is replaced into equation (4.4.15), after some calculations we get

$$(1-x^2)G_n''(x) - 2xG_n'(x) + \left((n+a+b)(n+a+b+1) - \frac{b^2}{1-x^2} + \frac{a((-1)^n - a)}{x^2}\right)G_n(x) = 0. \qquad (4.15.1)$$

By comparing (4.15.1) with (3.13) one can conclude that the sequence

$$Q_n(x;b^2, -\frac{2}{x^2}) = G_n(x;1,b) \quad ; \quad b > -1, \qquad (4.15.2)$$

is a general solution of equation (3.13) for $\mu_n = (n+b+1)(n+b+2)$, $\upsilon = b^2$ and $E(x) = -2/x^2$. Moreover, substituting these values into (3.14) and noting (4.14.3) yields (4.15.3)

$$\int_{-1}^{1} Q_n(x;b^2,-\frac{2}{x^2})Q_m(x;b^2,-\frac{2}{x^2})\,dx = \left(\prod_{i=1}^{n} \frac{(i+1-(-1)^i)(i+1-(-1)^i+2b)}{(2i+2b+1)(2i+2b+3)}\right)\frac{\sqrt{\pi}\,\Gamma(b+1)}{2\Gamma(b+5/2)}\delta_{n,m}$$

### 4.b. Second subclass: generalized Hermite polynomials (GHP)

GHP differential equation is a special case of the generic equation (4.7) for $\lambda = 2$ and consequently the equation (4.1) for $p = 0$, $q = 1$, $r = -2$, $s = 2u$, so that we have

$$x^2\Phi_n''(x) - 2x(x^2 - u)\Phi_n'(x) + \left(2nx^2 + ((-1)^n - 1)u\right)\Phi_n(x) = 0. \qquad (4.16)$$

According to (4.3) and (4.3.1), the explicit form of monic GHP is as (4.17)

$$\bar{S}_n\left(\begin{array}{cc} -2 & 2u \\ 0 & 1 \end{array} \bigg| x\right) = \prod_{i=0}^{[n/2]-1} \frac{2i + (-1)^{n+1} + 2 + 2u}{-2} \sum_{k=0}^{[n/2]} \binom{[n/2]}{k} \left(\prod_{i=0}^{[n/2]-(k+1)} \frac{-2}{2i + (-1)^{n+1} + 2 + 2u}\right) x^{n-2k}$$

that corresponds to the weight function

$$W\left(\begin{array}{cc} -2, & 2u \\ 0, & 1 \end{array} \bigg| x\right) = x^{2u} e^{-x^2}. \qquad (4.17.1)$$

Substituting (4.17) and (4.17.1) in the main theorem 1 when $\theta = \infty$ eventually yields (4.17.2)

$$\int_{-\infty}^{\infty} x^{2u} e^{-x^2} \bar{S}_n\left(\begin{array}{cc} -2 & 2u \\ 0 & 1 \end{array} \bigg| x\right) \bar{S}_m\left(\begin{array}{cc} -2 & 2u \\ 0 & 1 \end{array} \bigg| x\right) dx = \left((-1)^n \prod_{i=1}^{n} C_i\left(\begin{array}{cc} -2 & 2u \\ 0 & 1 \end{array}\right)\right)\Gamma(u+\frac{1}{2})\,\delta_{n,m}$$

Clearly (4.17.2) is valid when $u + 1/2 > 0$ and $(-1)^{2u} = 1$. To study GHP individually see [5,6,10].

**4.c. Third subclass: a finite class of symmetric orthogonal polynomials on** $(-\infty, \infty)$

Essentially there are two finite classes of symmetric orthogonal polynomials that are particular sub-cases of polynomials (4.3) and can be considered as two independent generalized Sturm-Lioville problems. To derive the first finite sub-class, first we define the initial vector

$$(a,b,c,d,\lambda) = (1,1,-2u-2v+2,-2u,2), \tag{4.18}$$

which gives the following symmetric differential equation (4.18.1)

$$x^2(x^2+1)\Phi_n''(x) - 2x((u+v-1)x^2+a)\Phi_n'(x) + \left(n(2u+2v-(n+1))x^2 + (1-(-1)^n)u\right)\Phi_n(x) = 0$$

The above differential equation is evidently a special case of the generic equation (4.7) and has an explicit polynomial solution as

$$\tag{4.18.2}$$

$$\bar{S}_n\left(\begin{array}{cc}-2u-2v+2, & -2u \\ 1, & 1\end{array}\middle| x\right) = \prod_{i=0}^{[n/2]-1} \frac{2i+(-1)^{n+1}+2-2u}{2i+2[n/2]+(-1)^{n+1}+2-2u-2v}$$

$$\times \sum_{k=0}^{[n/2]} \binom{[n/2]}{k} \left(\prod_{i=0}^{[n/2]-(k+1)} \frac{2i+2[n/2]+(-1)^{n+1}+2-2u-2v}{2i+(-1)^{n+1}+2-2u}\right) x^{n-2k}.$$

But, according to theorem 1, if (4.18.1) is written in self-adjoint form, then the following term must vanish on $(-\infty, \infty)$, i.e.

$$[x^{-2u}(1+x^2)^{-v+1}(\Phi_n'(x)\Phi_m(x) - \Phi_m'(x)\Phi_n(x))]_{-\infty}^{\infty} = 0, \tag{4.18.3}$$

where in fact

$$x^{-2u}(1+x^2)^{-v} = W\left(\begin{array}{cc}-2u-2v+2, & -2u \\ 1 & 1\end{array}\middle| x\right), \tag{4.18.4}$$

denotes the corresponding weight function. On the other hand, since $\Phi_n(x)$ is a polynomial of degree $n$, we should have

$$Max \deg(\Phi_n'(x)\Phi_m(x) - \Phi_m'(x)\Phi_n(x)) = n+m-1. \tag{4.19}$$

Therefore, from (4.18.3) it follows that the condition

$$-2u + 2(-v+1) + n + m - 1 \le 0, \tag{4.19.1}$$

must be valid in order that one can use the consequences of theorem 1. In other words, one should select one of the two following options

$$\begin{cases}-2u+n+m-1\le 0 \\ -2v+2\le 0\end{cases}, \quad \begin{cases}-2v+n+m-1\le 0 \\ -2u+2\le 0\end{cases} \tag{4.19.2}$$

in order to make a symmetric sequence of orthogonal polynomials. This yields the following corollary.

**Corollary 1.** The finite polynomial set $\{S_n(p,q,r,s;x)\}_{n=0}^{n=N}$ for $p=1$, $q=1$, $r=-2u-2v+2$ and $s=-2u$ is orthogonal with respect to the weight function $x^{-2u}(1+x^2)^{-v}$ on $(-\infty,\infty)$ if and only if $v \geq 1$, $N \leq u+1/2$ (or conversely $u \geq 1$, $N \leq v+1/2$).

According to above corollary and using Favard theorem, these polynomials should satisfy the relation

(4.20)

$$\int_{-\infty}^{\infty} \frac{x^{-2u}}{(1+x^2)^v} \bar{S}_n\left(\begin{array}{cc} -2u-2v+2, & -2u \\ 1, & 1 \end{array}\bigg| x\right) \bar{S}_m\left(\begin{array}{cc} -2u-2v+2, & -2u \\ 1, & 1 \end{array}\bigg| x\right) dx =$$

$$\left((-1)^n \prod_{i=1}^{n} C_i\left(\begin{array}{cc} -2u-2v+2, & -2u \\ 1, & 1 \end{array}\right)\right) \frac{\Gamma(u+v-1/2)\Gamma(-u+1/2)}{\Gamma(u+v)} \delta_{n,m}$$

iff $v \geq 1$ and $N = Max\{m,n\} \leq u+1/2$ or $u \geq 1$ and $N = Max\{m,n\} \leq v+1/2$.

**4.d. Fourth subclass: another finite class of symmetric orthogonal polynomials on $(-\infty,\infty)$.**

Similar to previous case, we can introduce another finite sequence of symmetric orthogonal polynomials on the interval $(-\infty,\infty)$ having the following main properties, see [8] for more details.

Initial vector: (4.21)
$$(a,b,c,d,\lambda) = (1,0,-2u+2,2,2),$$

Differential equation: (4.21.1)

$$x^4 \Phi_n''(x) + 2x((1-u)x^2 + 1)\Phi_n'(x) - \left(n(n+1-2u)x^2 + 1 - (-1)^n\right)\Phi_n(x) = 0.$$

Explicit form of monic polynomials: (4.21.2)

$$\bar{S}_n\left(\begin{array}{cc} -2u+2 & 2 \\ 1 & 0 \end{array}\bigg| x\right) = \prod_{i=0}^{[n/2]-1} \frac{2}{2i+2[n/2]+(-1)^{n+1}+2-2u} \times$$

$$\sum_{k=0}^{[n/2]} \binom{[n/2]}{k} (\prod_{i=0}^{[n/2]-(k+1)} \frac{2i+2[n/2]+(-1)^{n+1}+2-2u}{2}) x^{n-2k}.$$

Weight function:
$$W\left(\begin{array}{cc} -2u+2 & 2 \\ 1 & 0 \end{array}\bigg| x\right) = x^{-2u} \exp(-\frac{1}{x^2}).$$ (4.21.3)

Orthogonality relation:

$$\int_{-\infty}^{\infty} x^{-2u} \exp(-\frac{1}{x^2}) \bar{S}_n \left( \begin{array}{cc} -2u+2 & 2 \\ 1 & 0 \end{array} \Big| x \right) \bar{S}_m \left( \begin{array}{cc} -2u+2 & 2 \\ 1 & 0 \end{array} \Big| x \right) dx$$
$$= \left( (-1)^n \prod_{i=1}^{n} C_i \left( \begin{array}{cc} -2u+2 & 2 \\ 1 & 0 \end{array} \right) \right) \Gamma(u - \frac{1}{2}) \, \delta_{n,m}$$
(4.21.4)

if and only if $N = Max\{m, n\} \leq u - 1/2$.

As it was observed, each of polynomial classes of sections (4.a) to (4.d) were particular cases of a symmetric generalized Sturm-Liouville problem. However since all weight functions of presented polynomials were even functions, the condition $(-1)^{2u} = 1$ is always necessary. In other words, we can also consider the weights $|x|^{2u}(1-x^2)^v$, $|x|^{2u} e^{-x^2}$, $|x|^{-2u}(1+x^2)^{-v}$ and $|x|^{-2u} e^{-1/x^2}$ for four studied classes of symmetric orthogonal polynomials respectively.

**Remark 2.** Although some special functions such as Bessel [3], Fourier trigonometric sequences and so on are symmetric and satisfy a differential equation whose coefficients are alternatively even and odd, it is anyway important to note that they do not satisfy the conditions of theorem 1. For instance, if we choose

$$A(x) = x^2 \quad ; \quad B(x) = x \quad ; \quad C(x) = 1,$$
$$D(x) = x^2 \quad ; \quad E(x) = 0 \quad ; \quad \lambda_n = -n^2,$$
(4.22)

in the main equation (2.1), we get to the Bessel differential equation

$$x^2 \Phi_n''(x) + x \Phi_n'(x) + (x^2 - n^2) \Phi_n(x) = 0,$$
(4.22.1)

with a symmetric solution as

$$(-1)^n J_n(-x) = J_n(x) = \sum_{k=0}^{\infty} \frac{(-1)^k}{k!(n+k)!} (\frac{x}{2})^{n+2k} \quad , \quad n \in Z^+,$$
(4.22.2)

whereas the orthogonality interval of Bessel functions is not symmetric (i.e. $[0,1]$). Hence, the theorem 1 cannot be applied for $J_n(x)$, unless there exists a specific even function $E(x)$ for equation (2.1) such that the corresponding solution has infinity zeros (see e.g. [3] for more details), exactly similar to the case of usual Bessel functions of order $n$. This subject can be investigated as an open problem.

**5. Conclusion.**

In this paper, we showed that the usual Sturm-Liouville problem

$$\frac{d}{dx}\left(k(x)\frac{dy}{dx}\right) + (\lambda_n \rho(x) - q(x))y = 0 \quad ; \quad k(x) > 0, \; \rho(x) > 0,$$

under the boundary conditions

$$\begin{cases} \alpha_1 y(a) + \beta_1 y'(a) = 0 \\ \alpha_2 y(b) + \beta_2 y'(b) = 0 \end{cases} \quad ; \quad a < x < b$$

can be extended to the following problem

$$\frac{d}{dx}\left(k(x)\frac{dy}{dx}\right) + (\lambda_n \rho(x) - q(x) + \frac{1-(-1)^n}{2}E(x))y = 0,$$

with the boundary condition: $\alpha_3 y(\theta) + \beta_3 y'(\theta) = 0$ ; $-\theta < x < \theta$, provided that the solution of latter differential equation has symmetry property, i.e. $y_n(-x) = (-1)^n y_n(x)$.
Then we presented two basic examples to show the advantage of extending Sturm-Liouville problem for symmetric functions.
Finally let us add that the extension of Sturm-Liouville theorem for discrete variables is also possible (in preparation).

**References.**